\documentclass[preprint,3p,12pt,pdf]{elsarticle}

\usepackage{lineno,hyperref}
\modulolinenumbers[5]

\usepackage{hyperref}

\usepackage{epstopdf}

\usepackage{tikz}
\usetikzlibrary{arrows}

\usepackage{pst-node}
\usepackage{tikz-cd}

\usepackage[shellescape]{gmp}

\usepackage{mathrsfs}
\usepackage{amssymb}
\usepackage{amsfonts}
\usepackage{latexsym}
\usepackage{mathtools}
\usepackage{xcolor}
\usepackage{wrapfig}
\usepackage{floatflt}

\usepackage{mathtools}
\usepackage{extarrows}

\usepackage{graphicx}

\usepackage{subcaption}

\def\lb{\label}

\newcommand{\er}[1]{\textrm{(\ref{#1})}}

\newtheorem{theorem}{\bf Theorem}[section]

\def\a{\alpha}         
\def\b{\beta}          
         
\def\G{\Gamma}         
\def\d{\delta}

           \def\mJ{{\mathscr J}}
        
\def\k{\kappa}

\def\m{\mu}            
\def\n{\nu}

\def\ve{\varepsilon}           

    \def\R{{\mathbb R}}   \def\C{{\mathbb C}}    
    \def\N{{\mathbb N}}   

\def\lt{\biggl}                  \def\rt{\biggr}
               \def\wt{\widetilde}
\def\no{\noindent}

\let\ge\geqslant                 \let\le\leqslant

\def\iy{\infty}
\def\sm{\setminus}               
                 \def\ts{\times}
\def\pa{\partial}

\def\el2{\ell^{\,2}}             \def\1{1\!\!1}

\def\arg{\mathop{\mathrm{arg}}\nolimits}

\def\diag{\mathop{\mathrm{diag}}\nolimits}

\def\Im{\mathop{\mathrm{Im}}\nolimits}

\def\Ker{\mathop{\mathrm{Ker}}\nolimits}

\def\Re{\mathop{\mathrm{Re}}\nolimits}

\let\ge\geqslant
\let\le\leqslant

\newcommand{\ca}{\begin{cases}}
	\newcommand{\ac}{\end{cases}}
\newcommand{\ma}{\begin{pmatrix}}
	\newcommand{\am}{\end{pmatrix}}
\renewcommand{\[}{\begin{equation}}
	\renewcommand{\]}{\end{equation}}
\def\eq{\begin{equation}}
	\def\qe{\end{equation}}

\begin{document}
	
	\begin{frontmatter}

		\title{Complete left tail asymptotic for supercritical multitype branching processes}

		\date{\today}

		\author
		{Anton A. Kutsenko}
	
\address{University of Hamburg, MIN Faculty, Department of Mathematics, 20146 Hamburg, Germany; email: akucenko@gmail.com}
	
	\begin{abstract}
		We derive a complete left-tail asymptotic series for the density of the {\it martingale limit} of a supercritical multitype Galton-Watson process in the Schr\"oder case. We show that the series converges everywhere, not only for small arguments. This is the first complete result regarding the left tails of multitype branching processes. A good, quickly computed approximation for the density will also be derived from the series.
	\end{abstract}


	
\end{frontmatter}


{\section{Introduction}\lb{sec0}}

It is good to start with the classical single-type branching processes 
$$
 X_{t+1}=\sum_{j=1}^{X_t}\xi_{j,t},\ \ \ X_0=1,\ \ \ t\in\N\cup\{0\},
$$
where all $\xi_{j,t}$ are independent and identically distributed natural number-valued random variables.  We will only consider the supercritical Schr\"oder case, where the minimum family size is exactly one. In \cite{H}, some limit theorems and asymptotics are discussed. In particular, the {\it martingale limit} is defined by 
$$
 W:=\lim_{t\to+\iy}E^{-t}X_t,
$$ 
where $E$ is the expectation. The density $d(x)$ of $W$ is our primary object of discussion. For the small argument of the density, the first asymptotic term was analysed in \cite{D1}, where precise estimates were obtained. Then in \cite{BB1}, the authors found an explicit form of the first asymptotic term. After that, in \cite{K24}, the complete left tail asymptotic was derived. All the asymptotic terms were expressed explicitly through the so-called one-periodic Karlin-McGregor function introduced in \cite{KM1} and \cite{KM2}. Moreover, it was shown that the series converges everywhere, not only for small arguments of the density. To sum up the above, the main result of \cite{K24} is the following 
\[\lb{01}
 d(x)=\sum_{m=1}^{+\iy}x^{m\n-1}V_{m}(\log_E x),\ \ \ x>0,
\]
where $V_{m}(x)$ are one-periodic analytic functions, and $\n>0$ is some explicit parameter. All components of RHS in \er{01} are computable using efficient numerical procedures.

The $N$-type Galton-Watson process is defined as
$$
 {\bf X}_{n,t+1}=\sum_{i=1}^N\sum_{j=1}^{X_{ni,t}}\pmb{\xi}_{ij,t},\ \ \ {\bf X}_{n,0}=(\d_{jn})_{j=1}^N,\ \ \ n=1,...,N,\ \ \ t\in\N\cup\{0\},
$$
where all $\pmb{\xi}_{ij,t}$  are independent, and identically distributed for the same $i$. In this case, the variables are vectors of $N$ components representing different types of individuals. Each individual can produce individuals of different types simultaneously.
Following the Kesten-Stigum theorem, see \cite{KS} and \cite{AN}, one can define the martingale limits
$$
 {\bf W}_n=W_n{\bf a}:=\lim_{t\to+\iy}E^{-t}{\bf X}_{n,t}
$$
distributed along the half-line directed by the left eigenvector ${\bf a}$ of the mean matrix, corresponding to the maximal eigenvalue $E>0$. The densities ${\bf d}(x)=(d_n(x))_{n=1}^N$ of $(W_n)_{n=1}^N$ is our main object to study. The first asymptotic term of the logarithm of the density was discovered in \cite{J}. In this article, we derive the complete series
\[\lb{02}
{\bf d}(x)=\sum_{{\bf m}\in(\N\cup\{0\})^N\sm{\bf 0}}x^{{\bf m}\cdot\pmb{\n}-1}{\bf V}_{\bf m}(\log_E x),\ \ \ x>0,
\]
where ${\bf V}_{\bf m}(x)$ are one-periodic analytic vector-valued functions, and $\pmb{\n}$ is a vector with positive real parts. The corresponding main result with a rate of convergence of the series is formulated in Theorem \ref{T1} below. Again, all components of RHS in \er{02} are computable using efficient numerical procedures.  Both densities in \er{01} and \er{02} depend on the one-dimensional parameter $x>0$. The significant difference between the formulas is the powers of $x$. While powers of $x$ in \er{01} are regular, in \er{02}, they can be quasi-chaotic if some components of $\pmb{\n}$ are algebraically independent.  

The proof of \er{02} consists of several steps. It will be based on the analog of the inversion formula of Poincar\'e function \er{010} similar to that of single-type branching processes given in \cite{D1}. This section is discussed briefly because it primarily repeats the ideas presented in \cite{KS}, \cite{AN}, \cite{D1}, and to a lesser extent, \cite{J1}. The main part is the analysis of \er{010}. We try to follow the ideas of \cite{K24}, but some of the new methods should be developed to construct analogs of Karlin-McGregor functions. The use of a solution of the multidimensional Schr\"oder functional equation plays here one of the key roles. In comparison with the one-dimensional case, the existence of the solution requires more assumptions. In this regard, we use the results of \cite{Z}. Other sources of interesting analytical ideas useful in the analysis include \cite{WW}, \cite{FO}, \cite{M}, and \cite{K68}.

Usually, the oscillations of $V_m(x)$ in \er{01} are very small. They are ``almost'' constant functions.  Regarding functions similar to $V_1(x)$, such a phenomenon is known in physics, see, e.g., \cite{DIL}, \cite{CG}, and \cite{DMZ}. Since all $V_m(x)$ are dependent on $V_1(x)$, see \cite{K24}, the same phenomenon is expected for all of them. The situation with ${\bf V}_{\bf m}(x)$ seems to be similar.  Thus, it is natural to provide the following approximation of the density
\[\lb{03}
 \tilde{\bf d}(x)=\sum_{{\bf m}\in(\N\cup\{0\})^N\sm{\bf 0}}\frac{x^{{\bf m}\cdot\pmb{\n}-1}e^{{\bf m}\cdot\ln{\bf K}(1)}}{\G({\bf m}\cdot\pmb{\n})}\pmb{\psi}_{\bf m},\ \ \ x>0,
\] 
which is easy to compute, because we need to calculate only small part of ${\bf V}_{\bf m}(\log_E x)$ related to ${\bf V}_{\bf m}({\bf 0})$, see \er{028} and \er{040} below. In the Examples section, we compare \er{03} with \er{010}. Some perspectives are discussed in Section \ref{sec3}

{\section{Main results}\lb{sec1}}

As already mentioned above, the first part of this section consists of a brief discussion of \er{010}.  We will mainly follow the ideas of \cite{KS}, \cite{AN}, \cite{D1}, and a little bit of \cite{J1}. In this part, it is implicitly assumed that the conditions of the Kesten-Stigum theorem are satisfied, one of the most important of which is that some power of the expectation matrix \er{007} is positive-definite. However, our main result Theorem \ref{T1} is formulated for ${\bf d}(x)$ given in \er{010}. In this case, we can omit the conditions of the Kesten-Stigum theorem, working with \er{010} directly. Nevertheless, some new conditions {\bf (A)}-{\bf (D)} should be satisfied.

Let 
\[\lb{000}
 {\bf X}_{t}=(X_{1,t},...,X_{N,t})^T
\] 
be a multitype Galton-Watson branching process
\[\lb{001}
 {\bf X}_{t+1}=\sum_{i=1}^N\sum_{j=1}^{X_{i,t}}\pmb{\xi}_{ij,t},\ \ \ t\in\N\cup\{0\},
\]
where all $\pmb{\xi}_{ij,t}$  are independent. The symbol ${ }^T$ means standard transpose - for uniformity, I will mostly consider vector columns. The probability generating function (PGF) of  $\pmb{\xi}_{ij,t}$ is
\[\lb{002}
 P_i({\bf z})=\sum_{{\bf k}\in(\N\cup\{0\})^N}p_i({\bf k})z_1^{k_1}...z_N^{k_N},\ \ \ {\bf z}=(z_1,...,z_N)^T,
\]
where $p_i({\bf k})$ is the probability that an individual of type $i$ has the offspring vector ${\bf k}=(k_1,...,k_N)^T$. We will consider only the so-called Schr\"oder case
\[\lb{003} 
 p_1({\bf 0})=...=p_N({\bf 0})=0,\ \ \ {\bf 0}=(0,...,0)^T,
\]
or, in other words, ${\bf P}({\bf 0})={\bf 0}$, where
\[\lb{004}
 {\bf P}({\bf z})=(P_1({\bf z}),...,P_N({\bf z}))^T.
\]
The following assumptions exclude trivial cases
\[\lb{005}
 p_i({\bf k})<1,\ \ \ \forall i,{\bf k};\ \ \ all\ P_i\ are\ non-linear.
\]
For simplicity, we assume also that ${\bf P}({\bf z})$ is analytic for $|{\bf z}|:=\max_i|z_i|<1+c$ and some $c>0$. As a PGF, ${\bf P}({\bf z})$ is analytic for $|{\bf z}|<1$, but we have required a slightly more extended domain of its definition. For the first reading, it would be better to keep in mind that ${\bf P}({\bf z})$ is entire that covers most of the practical applications.

{\bf (A)} The last assumption plus \er{003} and \er{005} form our first technical assumption.

\no Specifying the initial conditions
\[\lb{006}
 {\bf X}_0={\bf e}_n:=(\d_{jn})_{j=1}^N,
\]
we define $N$ basic processes ${\bf X}_{n,t}$, $n=1,...,N$. To define the {\it martingale limit}, we need the matrix of expectations
\[\lb{007}
 {\bf E}=\lt(\frac{\pa P_i}{\pa z_j}({\bf 1})\rt)_{i,j=1}^N,\ \ \ {\bf 1}:=(1,...,1)^T.
\] 
There is a second assumption:

{\bf (B)} All the entries of ${\bf E}$ are positive. 

\no Then, by Perron-Frobenius (P-F) theorem, the maximal eigenvalue $E$ of ${\bf E}$ is positive and single. Since ${\bf P}({\bf z})$ is a non-trivial, by {\bf (A)}, probability generating vector-function (with non-negative Taylor coefficients) satisfying ${\bf P}({\bf 0})={\bf 0}$ and ${\bf P}({\bf 1})={\bf 1}$, it is easy to see that $E>1$. Again, by P-F theorem, all the entries of the left ${\bf a}$ and right ${\bf b}$ eigenvectors corresponding to $E$ are positive. We choose them such that ${\bf a}\cdot{\bf b}=1$. In this article, $\cdot$ denotes the standard sum of vector components products, without conjugation in the complex case. The projector on the corresponding eigenspace is 
\[\lb{008}
 {\bf L}:=\lim_{t\to+\iy}E^{-t}{\bf E}^t={\bf b}{\bf a}^T.
\]
Generally speaking, {\bf (B)} can be weakened, since we only need that $E>1$ is single and all the entries of the corresponding eigenvectors are positive. 

The {\it martingale limits} are defined by
\[\lb{009} 
 {\bf W}_n:=\lim_{t\to+\iy}E^{-t}{\bf X}_{n,t},\ \ \ n=1,...,N.
\]
All the {\it martingale limits} are distributed along the half-line $\{x{\bf a}:\ x>0\}$ and have the densities $d_n(x)$. They can be computed by the formula
\[\lb{010}
 {\bf d}(x):=(d_n(x))_{n=1}^N=\frac1{2\pi}\int_{-\iy}^{+\iy}\pmb{\Pi}(\mathbf{i}y)e^{\mathbf{i}yx}dy=\frac1{2\pi\mathbf{i}}\int_{\mathbf{i}\R}\pmb{\Pi}(z)e^{zx}dz,
\]
where $\pmb{\Pi}(z)$ is given by
\[\lb{011}
\pmb{\Pi}(z):=\lim_{t\to+\iy}\underbrace{{\bf P}\circ...\circ {\bf P}}_{t}({\bf 1}-{z}{E^{-t}}{\bf b}),
\]
it satisfies
\[\lb{012}
{\bf P}(\pmb{\Pi}(z))=\pmb{\Pi}(Ez),\ \ \ \pmb{\Pi}(0)={\bf 1},\ \ \ \pmb{\Pi}'(0)=-{\bf b}.
\]
Formulas \er{010}-\er{012} are analogs of those for single-type Galton-Watson processes derived in, e.g., \cite{D1}. At this stage, we will perceive equality \er{010} formally. The existence and convergence of the integral will be seen below, after we accept condition {\bf (C)}. At the beginning, let us discuss standard preliminary properties of the Poincar\'e-type mapping $\pmb{\Pi}(z)$. We mostly follow the ideas of \cite{K24} and \cite{M}. The domain of definition of $\pmb{\Pi}(z)$ depends on the properties of the zero-component 
\[\lb{013}
 \mJ_0:=\{{\bf z}\in\C^N:\ \underbrace{{\bf P}\circ...\circ {\bf P}}_{t}({\bf z})\to{\bf 0}\ for\ t\to+\iy\}
\]
of the filled Julia set for the mapping ${\bf P}({\bf z})$. As a non-trivial PGF with a zero free-term, ${\bf P}({\bf z})$ is contracting ($|{\bf P}({\bf z})|<|{\bf z}|$) inside the open unit ball $\{{\bf z}:\ |{\bf z}|<1\}$, and it has here a single attracting point ${\bf z}={\bf 0}$, where ${\bf P}({\bf 0})={\bf 0}$. In other words, self-iterations of ${\bf P}$ tend to zero inside the open unit ball. Hence, $\mJ_0$ contains the open unit ball. Thus, by \er{011} and standard arguments of holomorphic dynamics, $\pmb{\Pi}(z)$ is defined and analytic at least for $\Re z>0$, because $E>1$ and all the components of ${\bf b}$ are positive. Moreover, $|\pmb{\Pi}(z)|\to0$ if $z\to\iy$ and $|\arg z|<\pi/2-\ve$ for any fixed $\ve>0$. A small note: one way to prove analyticity of $\pmb{\Pi}(z)$ can be based on rapidly converging algorithms, such as \er{045}-\er{047}. Another way can be based on the functional equation \er{012} with the estimation of Taylor coefficients of $\pmb{\Pi}(z)$.  

It is convenient to repeat the above in terms of another set. By \er{011}, the properties of $\pmb{\Pi}(z)$ depend on the geometry of the modified filled Julia set
\[\lb{014}
 \mJ_1:=\{z\in\C:\ \underbrace{{\bf P}\circ...\circ {\bf P}}_{t}({\bf 1}-{z}{\bf b})\to{\bf 0}\ for\ t\to+\iy\}
\]
near the point $z=0$. The critical angle
\[\lb{015}
 \theta^*=\lim_{R\to+0}\sup\{\theta:\ re^{\mathbf{i}\wt\theta}\in\mJ_1\ for\ all\ |\wt\theta|<\theta\ and\ 0<r<R\}
\]
determines the sector, where $\pmb{\Pi}(z)$ decays. Namely, by \er{011}, we have that $\pmb{\Pi}(z)$ is defined and analytic in the sector $\{z\in\C:\ |\arg z|<\theta^*\}$ and $\pmb{\Pi}(z)\to{\bf 0}$ for $z\to\iy$ lying inside the sector and separated from its boundary by $|\arg z|<\theta^*-\ve$, for any fixed $\ve>0$. The exact power-law rate of decay of $\pmb{\Pi}(z)$ can be seen from \er{036} and \er{037} below. Due to {\bf (A)}, $\pmb{\Pi}(z)$ is also analytic in a neighborhood of $0$. The set $\mJ_1$ usually has a fractal boundary and can be very complex. However, since $\mJ_1$ is a special ``linear'' transformation of $\mJ_0$, it always contains some disc touching the origin and having a real center, which gives $\theta^*\ge\pi/2$. We require a little more:

{\bf (C)} We assume that the critical angle $\theta^*>\pi/2$.

\begin{figure}
	\centering
	\begin{subfigure}[b]{0.47\textwidth}
		\includegraphics[width=\textwidth]{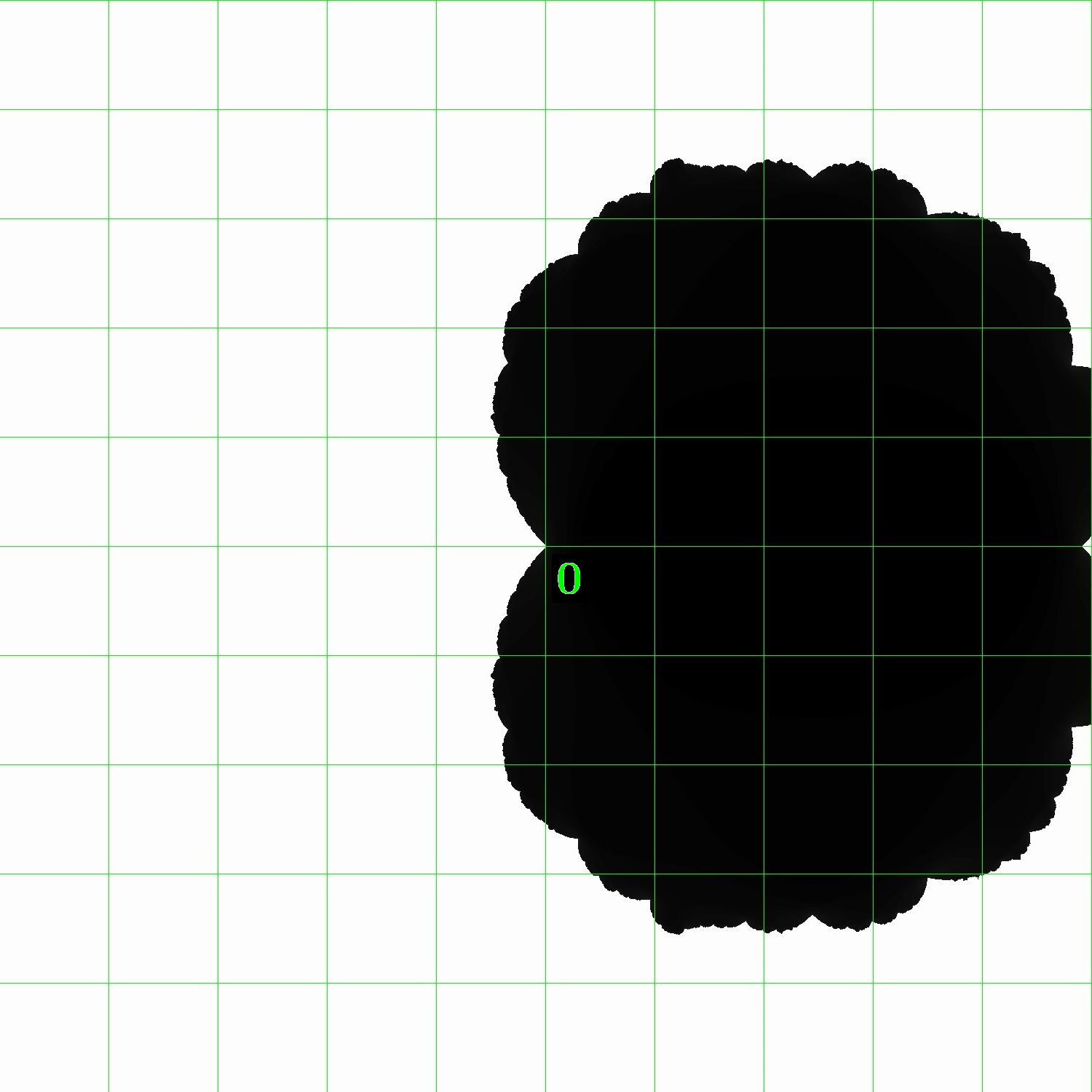}
		\caption{$\mJ_1$}
		\label{fig8a}
	\end{subfigure}
	~ 
	\begin{subfigure}[b]{0.47\textwidth}
		\includegraphics[width=\textwidth]{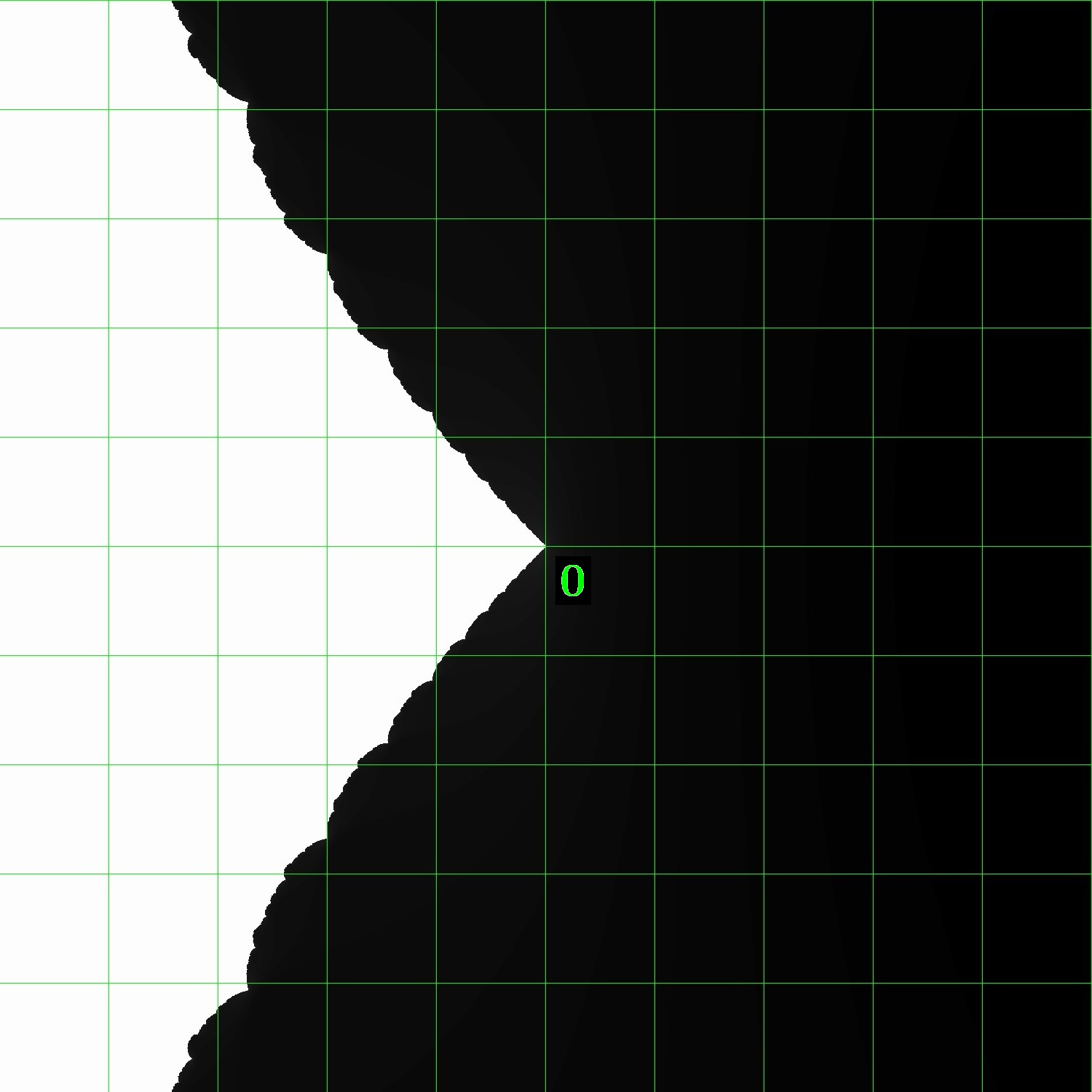}
		\caption{$\mJ_1$ zoomed in 10 times}
		\label{fig8b}
	\end{subfigure}
	\hfill
	\begin{subfigure}[b]{0.47\textwidth}
		\includegraphics[width=\textwidth]{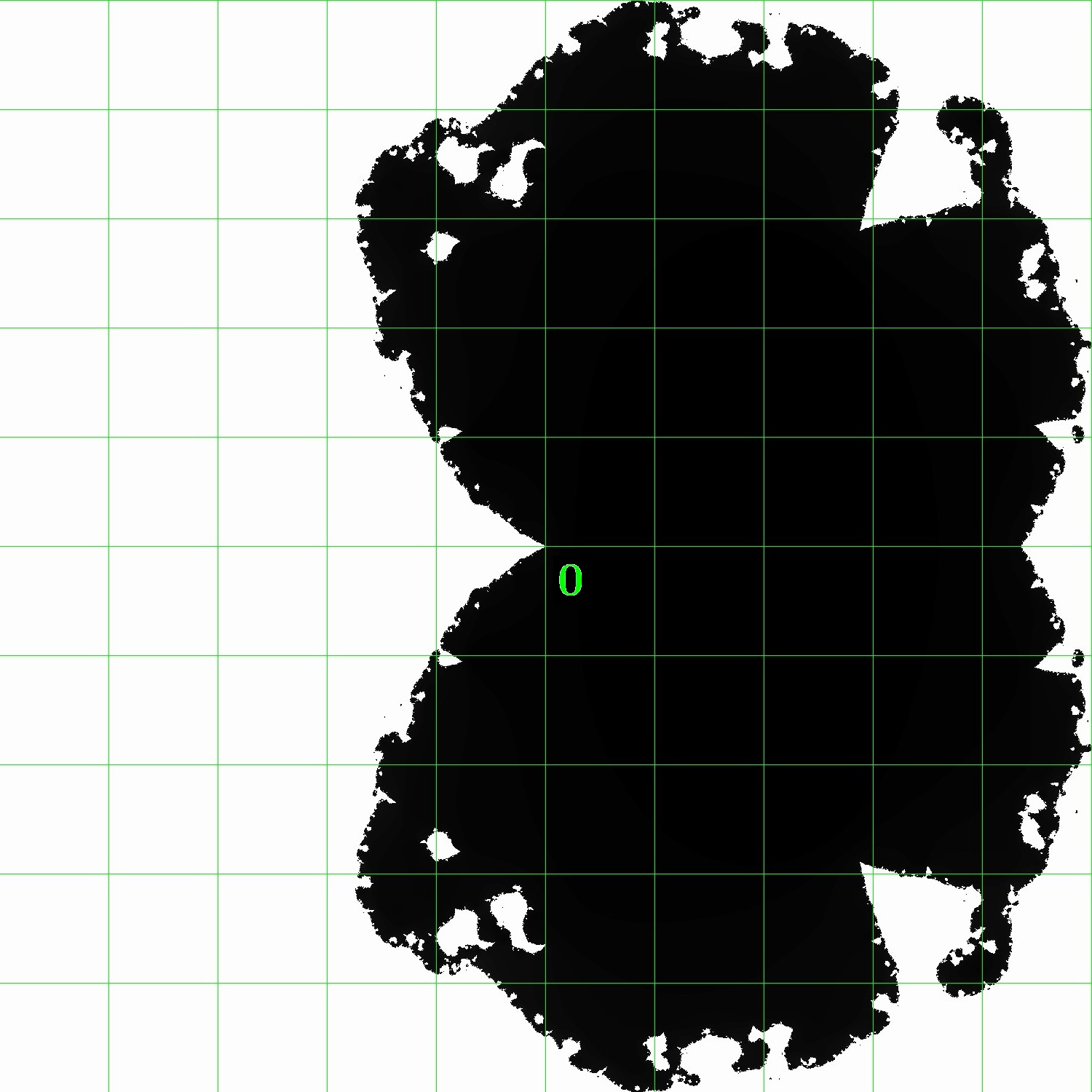}
		\caption{$\mJ_1$}
		\label{fig8c}
	\end{subfigure}
	~ 
	\begin{subfigure}[b]{0.47\textwidth}
		\includegraphics[width=\textwidth]{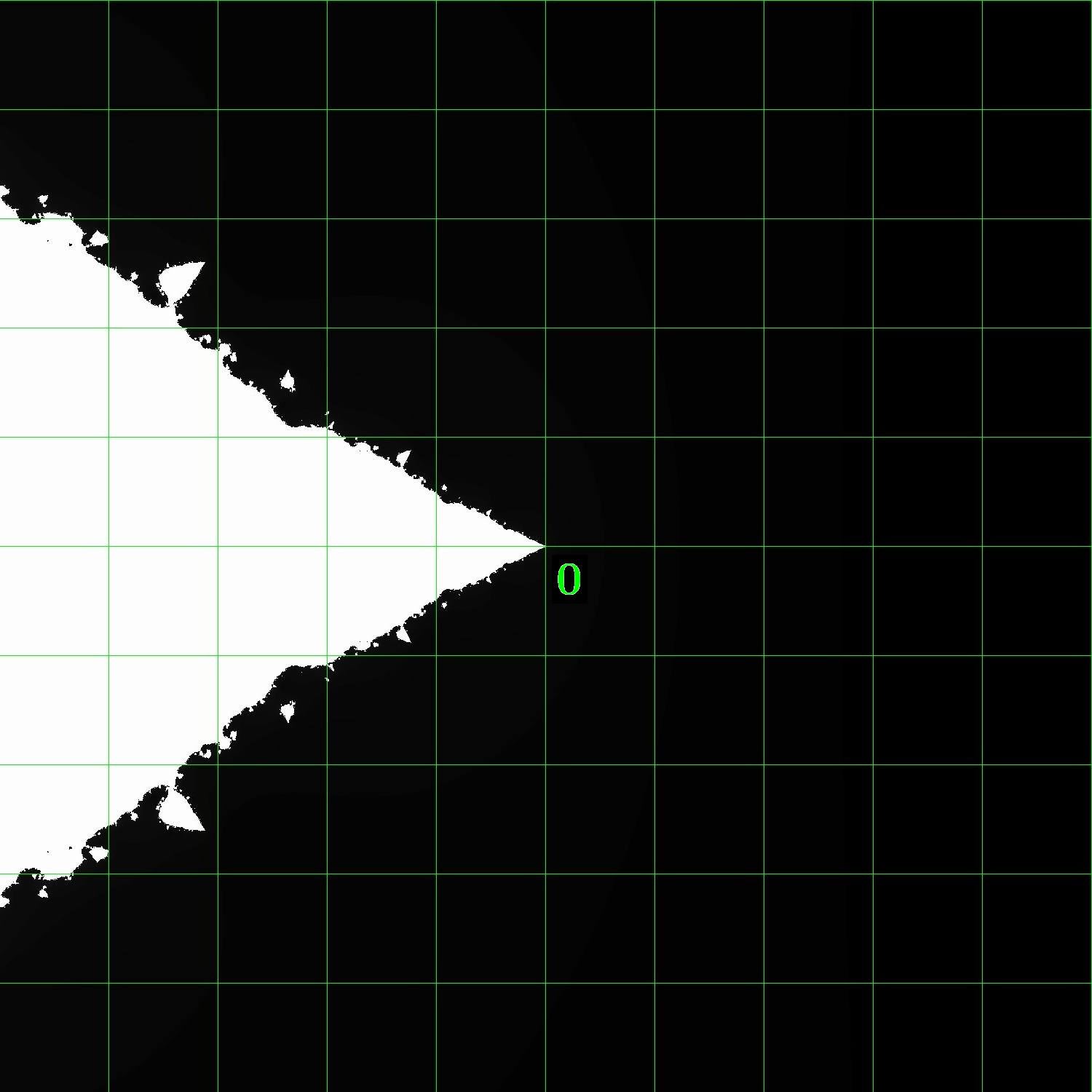}
		\caption{$\mJ_1$ zoomed in 10 times}
		\label{fig8d}
	\end{subfigure}
	\caption{The modified Julia set $\mJ_1$, see \er{014}. For (a) and (b), the PGF is ${\bf P}({\bf z})=((2z_1+z_2^2)/3,(z_2+z_1^2)/2)^T$; for (c) and (d), the PGF is ${\bf P}({\bf z})=((z_1+z_2+z_1^3)/(5-z_1z_2-z_1),0.3z_2+0.1z_1+0.3z_1z_2+0.1z_2^2+0.2z_1^3)^T$. It is seen that in both cases, the critical angle \er{015} is strictly greater than $\pi/2$.}\label{fig1}
\end{figure}

\no For all the examples we test, always $\theta^*>\pi/2$ strictly. The detailed analysis of the behavior of ${\bf P}({\bf z})$ near ${\bf z}={\bf 1}$ may avoid the condition {\bf (C)}, but this stays out of the main scope of the article. Two examples in Fig. \ref{fig1} show how simple or complex the form of $\mJ_1$ can be. As it is seen, the critical angle $\theta^*>\pi/2$ in both examples.

Due to {\bf (C)} and the arguments provided above {\bf (C)}, we have that $\pmb{\Pi}(z)$ is small for large $z$ in the right half-plane $\Re z\ge0$. Applying Cauchy's theorem to \er{010}, we can write
\[\lb{015a}
 {\bf d}(x)=\frac1{2\pi\mathbf{i}}\int_{\d+\mathbf{i}\R}\pmb{\Pi}(z)e^{zx}dz,\ \ \ \forall\d>0.
\]
This shift of the integration contour will prove very useful when expanding the density into a series below. For the further analysis of \er{010}, we need to define a few objects. Introduce
\[\lb{016}
{\bf M}=\lt(\frac{\pa P_i}{\pa z_j}({\bf 0})\rt)_{i,j=1}^N.
\]
In addition to the technical assumption {\bf (A)}, we require a more detailed item

{\bf (D)} The matrix ${\bf M}$ is diagonalizable and it satisfies
\[\lb{017}
 {\bf M}={\bf C}\diag(\m_1,...,\m_N){\bf C}^{-1},\ \ \ 0<|\m_N|^2<|\m_1|\le...\le|\m_N|.
\]
If {\bf (D)} is true then one may define
\[\lb{018}
 \pmb{\Phi}({\bf z}):=\lim_{t\to+\iy}{\bf M}^{-t}\underbrace{{\bf P}\circ...\circ {\bf P}}_{t}({\bf z}),
\]
which, by \cite{Z}, is analytic in some neighborhood of ${\bf 0}$. Following the definition \er{013}, it can be extended analytically to ${\bf z}\in\mJ_0$. By \er{016}, the functional equation for $\pmb{\Phi}({\bf z})$ is
\[\lb{019}
 \pmb{\Phi}({\bf P}({\bf z}))={\bf M}\pmb{\Phi}({\bf z}),\ \ \ \pmb{\Phi}({\bf 0})={\bf 0},\ \ \ \pmb{\Phi}({\bf z})\sim{\bf z}\ for\ {\bf z}\to{\bf 0}.
\]
Two functional equations \er{012} and \er{017}, along with the diagonalization \er{015} allow us to define the quasi-multiplicatively periodic function
\[\lb{020}
 {\bf K}(z)=(K_i(z))_{i=1}^N:={\bf C}^{-1}\pmb{\Phi}(\pmb{\Pi}(z)),
\]
each component of which satisfies
\[\lb{021}
 K_i(Ez)=\m_iK_i(z).
\]
Due to the definitions of $\pmb{\Phi}({\bf z})$, $\pmb{\Pi}(z)$, ${\bf K}(z)$, see \er{018}, \er{011}, \er{020}, and reasoning below \er{015}, ${\bf K}(z)$ is defined at least in the sector $\{z:\ |\arg z|<\theta^*\}$, where it is analytic. By analogy with periodic functions, quasi-multiplicatively periodic functions can be expanded into modified Fourier series
\[\lb{022}
 K_i(z)=\sum_{n=-\iy}^{+\iy}\k_{in}z^{\frac{\ln\m_i+2\pi\mathbf{i}n}{\ln E}},
\]
where
\[\lb{023}
 \k_{in}=\int_{0+\mathbf{i}y}^{1+\mathbf{i}y}\m_i^{-z}K_i(E^z)e^{-2\pi\mathbf{i}nz}dz,\ \ \ \forall y:\ |y|<\frac{\theta^*}{\ln E}.
\]
The reason of \er{022} and \er{023} is that $\m_i^{-z}K_i(E^z)$ is one-periodic in the strip $\{|\Im z|<\theta^*/\ln E\}$. One can easily obtain the following estimates by taking appropriate $y$ in \er{023}. For any $\ve>0$, there is $C_{\ve}$ such that
\[\lb{024}
 |\k_{in}|\le C_{\ve}e^{\frac{-2\pi |n|(\theta^*-\ve)}{\ln E}},\ \ \ \forall i,n.
\]
In \er{024}, we choose as $C_{\ve}$ the maximum of the modulus of integrands in \er{023} (without $e^{-2\pi\mathbf{i}nz}$) on the corresponding segment of integration, where $y=(\theta^*-\ve)/\ln E$. For our purposes, we need analogs of expansions \er{022} and estimates \er{024} for various products of functions $K_i(z)$. We use standard notations
\[\lb{025}
 \pmb{\m}^{\pmb{\a}}:=\prod_{i=1}^N\m_i^{\a_i},\ \ \ \pmb{\m}=(\m_i)_{i=1}^N,\ \ \ \pmb{\a}=(\a_i)_{i=1}^N,\ \ \ \|\pmb{\a}\|:=\sum_{i=1}^N|\a_i|.
\]
Let us arrange all the products $\pmb{\m}^{\pmb{\a}}$ for $\pmb{\a}\in(\N\cup\{0\})^N\sm{\bf 0}$ in descending order of their modules
\[\lb{026}
 |\pmb{\m}^{\bf\hat 1}|\ge|\pmb{\m}^{\bf\hat 2}|\ge |\pmb{\m}^{\bf\hat 3}|\ge...,
\]
where 
${\bf\hat 1}$, ${\bf\hat 2}$, ${\bf\hat 3}$ and so on are the corresponding multiindices. Recall that, by \er{016} and \er{017}, all $|\m_i|<1$ because ${\bf P}({\bf z})$ is a non-trivial PGF. All the products of quasi-multiplicatively periodic functions are also quasi-multiplicatively periodic, and we can repeat arguments \er{021}-\er{024}, namely 
\[\lb{027}
 {\bf K}^{\bf\hat m}(z)=\sum_{n=-\iy}^{+\iy}\k_{{\bf\hat m}n}z^{\frac{\ln\pmb{\m}^{\bf\hat m}+2\pi\mathbf{i}n}{\ln E}},
\]
where
\[\lb{028}
 \k_{{\bf\hat m}n}=\int_{0+\mathbf{i}y}^{1+\mathbf{i}y}(\pmb{\m}^{\bf\hat m})^{-z}{\bf K}^{\bf\hat m}(E^z)e^{-2\pi\mathbf{i}nz}dz,\ \ \ \forall y:\ |y|<\frac{\theta^*}{\ln E},
\]
which satisfy
\[\lb{029}
 |\k_{{\bf\hat m}n}|\le C_{\ve}^{\|{\bf\hat m}\|}e^{\frac{-2\pi |n|(\theta^*-\ve)}{\ln E}},\ \ \ \forall {\bf\hat m},n,
\]
and any sufficiently small $\ve>0$. Due to \er{029}, if $\Re z>1$ then the series in \er{027} converges exponentially fast, since $\arg z<\pi/2$ and $\theta^*>\pi/2$. More precisely, applying in \er{027} standard estimates for the power functions in the complex plane and using \er{029}, we obtain
\begin{multline}\lb{030}
 |\k_{{\bf\hat m}n}z^{\frac{\ln\pmb{\m}^{\bf\hat m}+2\pi\mathbf{i}n}{\ln E}}|\le C_{\ve}^{\|{\bf\hat m}\|}e^{\frac{-2\pi |n|(\theta^*-\ve)}{\ln E}}\lt|e^{\frac{(\ln|z|+\mathbf{i}\arg z)(\sum_{i=1}^N(\ln|\m_i|+\mathbf{i}\arg\m_i){\hat m}_i+2\pi\mathbf{i}n)}{\ln E}}\rt|\\
=C_{\ve}^{\|{\bf\hat m}\|}e^{\frac{-2\pi |n|(\theta^*-\ve)}{\ln E}}e^{\frac{\ln|z|\sum_{i=1}^N{\hat m}_i\ln|\m_i|-2\pi n\arg z-\arg z\sum_{i=1}^N{\hat m}_i\arg\m_i}{\ln E}}\\
=(AC_{\ve}|z|^{\max_i\log_E|\m_i|})^{\|\hat{\bf m}\|}e^{\frac{-2\pi |n|(\theta^*-\frac{\pi}2-\ve)}{\ln E}},\ \ \ \Re z>1,
\end{multline}
with some constant $A>0$ depending on $\arg\m_i$. Recall that all $|\m_i|<1$. Let us fix some $\ve_0>0$ such that $\theta^*-{\pi}/2-\ve_0>0$, see the condition {\bf (C)} above. Then we can write \er{030} as
\[\lb{031}
 |\k_{{\bf\hat m}n}z^{\frac{\ln\pmb{\m}^{\bf\hat m}+2\pi\mathbf{i}n}{\ln E}}|\le(B|z|^{\b})^{\|\hat{\bf m}\|}\a^{|n|},\ \ \ \Re z>1,
\]
where the corresponding constants satisfy
\[\lb{032}
 B>0,\ \ \ \b<0,\ \ \ 0<\a<1.
\]
Substituting \er{031} into \er{027}, and using an explicit form of the sum of geometric series, we obtain
\[\lb{033}
 |{\bf K}^{\bf\hat m}(z)|\le (B|z|^{\b})^{\|{\bf\hat m}\|}\frac{1+\a}{1-\a},\ \ \ \Re z>1,
\]
which shows how small ${\bf K}^{\bf\hat m}(z)$ for large $z$.

The analytic function $\pmb{\Phi}({\bf z})\sim{\bf z}$ for small arguments, see \er{019}. Hence, it has an analytic inverse. Let us denote
\[\lb{034}
 \pmb{\Psi}({\bf z}):=\pmb{\Phi}^{-1}({\bf C}{\bf z})=\sum_{{\bf m}\in(\N\cup{\bf 0})^N\sm{\bf 0}}{\bf z}^{\bf m}\pmb{\psi}_{\bf m}.
\]
The series converges absolutely for $|{\bf z}|<r$,
with some $r>0$. The Cauchy estimates are
\[\lb{035}
 |\pmb{\psi}_{\bf m}|\le Dr^{-\|{\bf m}\|},
\]
with some constant $D>0$.

Estimates \er{033} show that ${\bf K}(z)$ is small for $\Re z>\d$ with large $\d>0$. Thus, ${\pmb \Psi}({\bf K}(z))$ is defined and analytic for such $z$. Using \er{020}, along with \er{034} and \er{027}, we obtain
\[\lb{036}
 \Pi(z)=\sum_{m=1}^{+\iy}\sum_{n=-\iy}^{+\iy}\k_{{\bf\hat m}n}z^{\frac{\ln\pmb{\m}^{\bf\hat m}+2\pi\mathbf{i}n}{\ln E}}\pmb{\psi}_{\bf\hat m},\ \ \ \Re z>\d,
\]
for all sufficiently large $\d>0$. Moreover, due to \er{031} and \er{033}, the series \er{036} converges absolutely. More precisely, let us fix some $m_0>0$ and estimate the reminder
\begin{multline}\lb{037}
 |\sum_{|{\bf\hat m}|\ge m_0}\sum_{n=-\iy}^{+\iy}\k_{{\bf\hat m}n}z^{\frac{\ln\pmb{\m}^{\bf\hat m}+2\pi\mathbf{i}n}{\ln E}}\pmb{\psi}_{\bf\hat m}|\le\sum_{|{\bf\hat m}|\ge m_0}\sum_{n=-\iy}^{+\iy}D(B|z|^{\b}r^{-1})^{\|\hat{\bf m}\|}\a^{|n|}\\
\le\frac{D(1+\a)(B|z|^{\b}r^{-1})^{Nm_0}}{(1-\a)(1-B|z|^{\b}r^{-1})^N}\le \frac{D(1+\a)(B|z|^{\b}r^{-1})^{Nm_0}}{(1-\a)(1-B\d^{\b}r^{-1})^N}=:R_{m_0}(z),\ \ \ \Re z>\d,
\end{multline}
for all sufficiently large $\d>0$. From \er{037}, we get
\[\lb{038}
 \frac1{2\pi\mathbf{i}}\int_{\d+\mathbf{i}\R}R_{m_0}(z)e^{zx}dz\to0,\ \ \ m_0\to+\iy,
\]
for any fixed $x>0$ and all sufficiently large $\d>0$, since $e^{zx}$ is bounded along lines parallel to the imaginary axis in the complex plane. Using \er{038}, along with the exponential fast convergence of \er{036} by the parameter $n$ for any fixed $m$, see \er{031}, one can change the order of summation and integration in the following identity, see \er{015a} and \er{036},
\begin{multline}\lb{039}
 {\bf d}(x)=\frac1{2\pi\mathbf{i}}\int_{\d+\mathbf{i}\R}\pmb{\Pi}(z)e^{zx}dz=\frac1{2\pi\mathbf{i}}\int_{\d+\mathbf{i}\R}\sum_{m=1}^{+\iy}\sum_{n=-\iy}^{+\iy}\k_{{\bf\hat m}n}z^{\frac{\ln\pmb{\m}^{\bf\hat m}+2\pi\mathbf{i}n}{\ln E}}e^{zx}\pmb{\psi}_{\bf\hat m}dz\\
=\sum_{m=1}^{+\iy}\sum_{n=-\iy}^{+\iy}\frac{\k_{{\bf\hat m}n}}{2\pi\mathbf{i}}\int_{\d+\mathbf{i}\R}z^{\frac{\ln\pmb{\m}^{\bf\hat m}+2\pi\mathbf{i}n}{\ln E}}e^{zx}dz\pmb{\psi}_{\bf\hat m}\\
=\sum_{m=1}^{+\iy}\sum_{n=-\iy}^{+\iy}\frac{\k_{{\bf\hat m}n}x^{-1-\frac{\ln\pmb{\m}^{\bf\hat m}+2\pi\mathbf{i}n}{\ln E}}}{2\pi\mathbf{i}}\int_{x\d+\mathbf{i}\R}z^{\frac{\ln\pmb{\m}^{\bf\hat m}+2\pi\mathbf{i}n}{\ln E}}e^{z}dz\pmb{\psi}_{\bf\hat m}\\
=\sum_{m=1}^{+\iy}\sum_{n=-\iy}^{+\iy}\frac{\k_{{\bf\hat m}n}x^{-1-\frac{\ln\pmb{\m}^{\bf\hat m}+2\pi\mathbf{i}n}{\ln E}}}{\G(-\frac{\ln\pmb{\m}^{\bf\hat m}+2\pi\mathbf{i}n}{\ln E})}\pmb{\psi}_{\bf\hat m},
\end{multline}
for all sufficiently large $\d>0$. In \er{039}, for the latest integral, we use the Hankel representation of the $\Gamma$-function, see Example 12.2.6 on p. 254 in \cite{WW}. Denoting
\[\lb{040}
 {\bf V}_{\bf\hat m}(x):=\sum_{n=-\iy}^{+\iy}\frac{\k_{{\bf\hat m}n}e^{-2\pi\mathbf{i}nx}}{\G(-\frac{\ln\pmb{\m}^{\bf\hat m}+2\pi\mathbf{i}n}{\ln E})}\pmb{\psi}_{\bf\hat m},
\]
we can formulate our main result
\begin{theorem}\lb{T1}
If conditions {\bf (A)}-{\bf (D)} are satisfied, then for ${\bf d}(x)$ given in \er{010}, the following expansion is true
\[\lb{041}
 {\bf d}(x)=\sum_{m=1}^{+\iy}x^{-1-{\log_E\pmb{\m}^{\bf\hat m}}}{\bf V}_{\bf\hat m}(\log_Ex),\ \ \ x>0,
\]
where one-periodic vector-valued functions ${\bf V}_{\bf\hat m}(x)$ are defined in \er{040}. Substituting \er{040} into \er{041}, we obtain the double series \er{039}, which converges exponentially fast by $n$ and super-exponentially fast by $m$.
\end{theorem}
The expansion \er{041} is proven above. To prove the rate of convergence discussed at the end of the theorem, we need to estimate the Gamma function. Using Stirling's approximation, see, e.g., \cite{WW},
$$
 \G(z)=\sqrt{\frac{2\pi}{z}}\lt(\frac{z}{e}\rt)^z\lt(1+O_{\ve}\lt(\frac1z\rt)\rt),\ \ \ |\arg z|<\pi-\ve,
$$
we obtain the following rough estimate:
$$
 if\ x>1\ and\ y\in\R\ then\ \ \frac1{|\G(x+\mathbf{i}y)|}\le D(1+|y|)e^{-sx\ln x+\frac{\pi|y| }2}\ \ for\ some\ D,s>0.
$$
Thus, for all sufficiently large $n$ and $m$, we obtain
\[\lb{042}
 \frac1{|\Gamma(-\frac{\ln\pmb{\m}^{\bf\hat m}+2\pi\mathbf{i}n}{\ln E})|}\le D_1|n|e^{-s_1\|{\bf\hat m}\|\ln\|{\bf\hat m}\|+\frac{\pi^2 |n|}{\ln E}}\ \ \ for\ some\ D_1,s_1>0.
\]
Combining \er{042} with \er{029} and \er{035}, we obtain
\[\lb{043}
 \lt|\frac{\k_{{\bf\hat m}n}x^{-1-{\log_E\pmb{\m}^{\bf\hat m}}}}{\G(-\frac{\ln\pmb{\m}^{\bf\hat m}+2\pi\mathbf{i}n}{\ln E})}\pmb{\psi}_{\bf\hat m}\rt|\le D_2e^{-s_2\|{\bf\hat m}\|\ln\|{\bf\hat m}\|-\d_2|n|}\ \ \ for\ some\ D_2,s_2,\d_2>0,
\]
for (any) fixed $x>0$ and all sufficiently large $n$ and $m$. Here, it is important that the critical angle $\theta^*>\pi/2$ in \er{029}. The estimate \er{043} finishes the proof of Theorem \ref{T1}.

All the components of \er{040} can be computed by using efficient numerical procedures. For the special functions \er{011} and \er{018} we can use algorithms similar to those developed in, e.g., \cite{K} and \cite{K24}. Then, the Fourier coefficients $\k_{{\bf\hat m}n}$, see \er{028}, can be computed with the help of Fast Fourier Transform. Due to \er{019} and \er{017}, the function $\pmb {\Psi}({\bf z})$, see \er{034}, satisfies the following functional equation
\[\lb{044}
 \pmb{\Psi}(\diag(\pmb{\m}){\bf z})={\bf P}(\pmb {\Psi}({\bf z})),\ \ \ \pmb{\Psi}({\bf z})\sim{\bf C}{\bf z}\ for\ {\bf z}\to{\bf 0}.
\]
This equation allows us to determine Taylor coefficients $\pmb{\psi}_{\bf m}$ step by step, substituting \er{034} into \er{044}. The ordering ${\bf m}\leftrightarrow{\bf\hat m}$ of multi-indices by \er{026} is not necessary, because of exponentially super-fast convergence of the series \er{041} by $m$. The ordering is used only to bring asymptotic properties to \er{041} for $x\to+0$, where each next term is asymptotically smaller than the previous one.   

Let us discuss how to compute $\pmb{\Pi}(z)$. It is convenient to extend the definition of $\pmb{\Pi}(z)$ by
\[\lb{045}
 \pmb{\Pi}_t({\bf z})=\underbrace{{\bf P}\circ...\circ {\bf P}}_{t}({\bf 1}-{E^{-t}}{\bf z}),\ \ \ \pmb{\Pi}({\bf z}):=\lim_{t\to+\iy}\pmb{\Pi}_t({\bf z}).
\]
Taylor expansion of PGF ${\bf P}({\bf z})$ at ${\bf z}={\bf 1}$ is
\[\lb{046}
 {\bf P}({\bf 1}-{\bf y})={\bf 1}-{\bf E}{\bf y}+{\bf P}_1({\bf y}),\ \ where\ \ |{\bf P}_1({\bf y})|\le P_1|{\bf y}|^2,
\]
for all sufficiently small ${\bf y}$ and some $P_1>0$. Thus, by \er{045} and \er{046}, we have
\[\lb{047}
 \pmb{\Pi}_{t+1}({\bf z})=\pmb{\Pi}_t(E^{-1}{\bf E}{\bf z}-E^{t}{\bf P}_1(E^{-t-1}{\bf z})),
\]
where the maximal by modulus eigenvalue of $E^{-1}{\bf E}$ is $1$, see {\bf (B)}, and $|E^{t}{\bf P}_1(E^{-t-1}{\bf z})|\le P_1E^{-t-2}|{\bf z}|$ is exponentially small for large $t$. Hence, \er{047} provides an efficient recurrence scheme for the computation of $\pmb{\Pi}({\bf z})$ - it shows an exponentially fast convergence, similar to that found in \cite{K}. Moreover, if ${\bf z}=z{\bf b}+{\bf w}$ with ${\bf w}\in\Ker{\bf L}$, see \er{008}, then $\pmb{\Pi}({\bf z})=\pmb{\Pi}(z{\bf b})$, since $E^{-1}{\bf E}$ is contracting for ${\bf w}\in\Ker{\bf L}$ by {\bf (B)}. At the end, $\pmb{\Pi}(z)=\pmb{\Pi}(z{\bf b})$ is our goal of computation.
 
The computation of $\pmb{\Phi}({\bf z})$ is based on computations of
\[\lb{048}
 \pmb{\Phi}_t({\bf z}):={\bf M}^{-t}\underbrace{{\bf P}\circ...\circ {\bf P}}_{t}({\bf z}).
\]
Taylor expansion of PGF ${\bf P}({\bf z})$ at ${\bf z}={\bf 0}$ is
\[\lb{049}
{\bf P}({\bf y})={\bf M}{\bf y}+{\bf Q}_1({\bf y}),\ \ where\ \ |{\bf Q}_1({\bf y})|\le Q_1|{\bf y}|^2,
\]
for all sufficiently small ${\bf y}$ and some $Q_1>0$. Then
\[\lb{050}
 \pmb{\Phi}_{t+1}({\bf z})=\pmb{\Phi}_{t}({\bf z})+{\bf M}^{-t-1}{\bf Q}_1({\bf M}^{t}\pmb{\Phi}_{t}({\bf z})),
\]
where
\[\lb{051}
 |{\bf M}^{-t-1}{\bf Q}_1({\bf M}^{t}\pmb{\Phi}_{t}({\bf z}))|\le\wt Q_1|\m_1|^{-t-1}|\m_N|^{2t}|\pmb{\Phi}_{t}({\bf z})|=\frac{\wt Q_1}{|\m_1|}\cdot\lt(\frac{|\m_N|^2}{|\m_1|}\rt)^t|\pmb{\Phi}_{t}({\bf z})|,
\]
for ${\bf z}\in\mJ_0$, see \er{013}. In \er{051}, the constant $\wt Q_1>0$ depends on $Q_1$ and the norms of ${\bf C}$ and ${\bf C}^{-1}$, see \er{017}. Due to {\bf (D)}, see again \er{017}, and \er{051}, the second term in RHS of \er{050} is exponentially small for large $t$. Thus, \er{050} provides an efficient recurrence scheme for the computation $\pmb{\Phi}({\bf z})=\lim_{t\to+\iy}\pmb{\Phi}_t({\bf z})$ - it shows an exponentially fast convergence.

{\section{Examples}\lb{sec1}}

For the calculations, Embarcadero Delphi Rad Studio Community Edition and the library NesLib.Multiprecision are used. This software provides a convenient environment for programming and well-functioning basic functions for high-precision computations. All the algorithms related to the article's subject, including efficient parallelization, are developed by the author (AK). We use 128-bit precision instead of the standard 64-bit double precision. The reason is the computation of Taylor coefficients $\pmb{\psi}_{\bf m}$ of $\pmb{\Psi}({\bf z})$, see \er{044}. However, it seems that the standard precision is enough for all the computations, because some techniques can stabilize the computation of $\pmb{\psi}_{\bf m}$. In this regard, the use of a quick 64-bit BLAS/LAPACK adaptation Dew.MtxVec instead of 128/256-bit NesLib.Multiprecision, as I did in previous articles, does not greatly reduce the accuracy. Nevertheless, if our goal is small oscillations in ${\bf V}_{\bf m}(x)$, then NesLib.Multiprecision is extremely helpful.

Let us consider the following two-dimensional PGF
\[\lb{200}
 {\bf P}({\bf z})=\ma pz_1+(1-p)z_2^2\\
                      qz_2+(1-q)z_1^2 \am,\ \ \ 0<p,q<1.
\]
The expectation matrix and its maximal eigenvalue are
\[\lb{201}
 {\bf E}=\ma 
              p & 2(1-p) \\
              2(1-q) & q
         \am,\ \ \ E=\frac{p+q+\sqrt{16(1-p)(1-q)+(p-q)^2}}{2}.
\]
The left and right eigenvectors of ${\bf E}$ are
\[\lb{202}
 {\bf a}=\frac1{\sqrt{16(1-p)(1-q)+(p-q)^2}}\ma \frac{E-q}{2(1-p)} \\ 1 \am,\ \ \ {\bf b}=\ma 2(1-p) \\ E-p \am,\ \ \ {\bf a}\cdot{\bf b}=1.
\]
The Taylor expansion of ${\bf P}({\bf z})$ at ${\bf z}={\bf 1}$ is
\[\lb{203}
 {\bf P}({\bf 1}-{\bf y})={\bf 1}-{\bf E}{\bf y}+{\bf P}_1({\bf y}),\ \ \ {\bf P}_1({\bf y})=\ma
                   (1-p)y_2^2 \\
                   (1-q)y_1^2
               \am.
\]
For the computation of $\pmb{\Pi}(z)=\pmb{\Pi}(z{\bf b})$, we apply the scheme based on \er{045} and \er{047}. The Taylor expansion of ${\bf P}({\bf z})$ at ${\bf z}={\bf 0}$ is
\[\lb{204}
 {\bf P}({\bf y})={\bf M}{\bf y}+{\bf Q}_1({\bf y}),\ \ \ {\bf M}=\diag(p,q),\ \ \ {\bf Q}_1({\bf y})={\bf P}_1({\bf y}).
\]
For the computation of $\pmb{\Phi}({\bf z})$, we use the scheme based on \er{048} and \er{050}. Both computations of $\pmb{\Pi}({\bf z})$ and $\pmb{\Phi}({\bf z})$ show exponentially fast convergence.

For the computations of $\pmb{\psi}_{\bf m}$, see \er{034}, the scheme based on the functional equation \er{044}, is applied. In our example, \er{044} takes the form
\[\lb{205}
 \ma{\Psi}_1(pz_1,qz_2)\\ {\Psi}_2(pz_1,qz_2)\am=\ma p\Psi_1(z_1,z_2)+(1-p)\Psi_2(z_1,z_2)^2\\
                      q\Psi_2(z_1,z_2)+(1-q)\Psi_1(z_1,z_2)^2 \am,\ \ \ \ma{\Psi}_1\\ {\Psi}_2\am\sim\ma z_1\\ z_2 \am,\ {\bf z}\to{\bf 0}.
\]
Substituting expansions
\[\lb{206}
 \Psi_1(z_1,z_2)=\sum_{{\bf m}\in(\N\cup{\bf 0})^N\sm{\bf 0}}\psi_{1,m_1m_2}z_1^{m_1}z_2^{m_2},\ \ \ \Psi_2(z_1,z_2)=\sum_{{\bf m}\in(\N\cup{\bf 0})^N\sm{\bf 0}}\psi_{2,m_1m_2}z_1^{m_1}z_2^{m_2}
\]
into \er{205}, we obtain
\[\lb{207}
 \psi_{1,{\bf m}}=\frac{1-p}{p^{m_1}q^{m_2}-p}\sum_{{\bf f}+{\bf g}={\bf m}}\psi_{2,{\bf f}}\psi_{2,{\bf g}},\ \ \ \psi_{2,{\bf m}}=\frac{1-q}{p^{m_1}q^{m_2}-p}\sum_{{\bf f}+{\bf g}={\bf m}}\psi_{1,{\bf f}}\psi_{1,{\bf g}}.
\]
Equation \er{207} along with the initial conditions
\[\lb{208}
 \psi_{1,10}=\psi_{2,01}=1,\ \ \ \psi_{1,01}=\psi_{2,10}=0,
\]
allows us to determine all the Taylor coefficients $\pmb{\psi}_{\bf m}$ step by step. They grow exponentially fast. To stabilize the computations, it is conveniet to compute $r^{\|{\bf m}\|}\pmb{\psi}_{\bf m}$, with some small $r$, instead of $\pmb{\psi}_{\bf m}$. To do this, it is enough to set $\psi_{1,10}=\psi_{2,01}=r$ in \er{208}. Then, in the computation of the final result \er{03}, one can apply the next technique
\[\lb{209}
 \frac{{\psi}_{i,{\bf m}}}{\G(-\log_E\pmb{\m}^{\bf m})}=\exp\lt(\ln r^{\|{\bf m}\|}{\psi}_{i,{\bf m}}-\ln\G(-{\bf m}\cdot\log_E\pmb{\m})-\|{\bf m}\|\ln r\rt),
\]
for $i=1,2$ in our case. The logarithm in \er{209} stabilizes also the computation of the Gamma function: for $\ln\G(z)$, there are many good numerical algorithms.

\begin{figure}
	\centering
	\begin{subfigure}[b]{0.75\textwidth}
		\includegraphics[width=\textwidth]{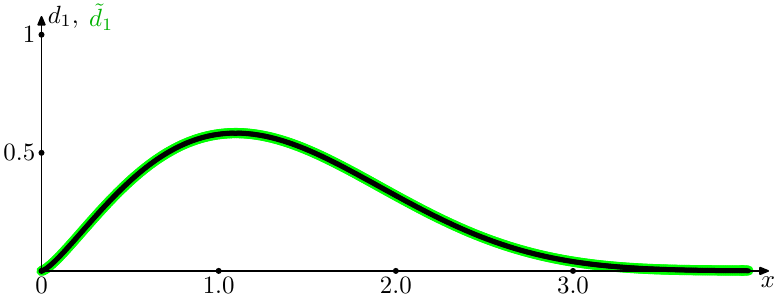}
		\caption{$p=1/3$, $q=1/2$}
		\label{fig1a}
	\end{subfigure}
	\hfill
	\begin{subfigure}[b]{0.75\textwidth}
		\includegraphics[width=\textwidth]{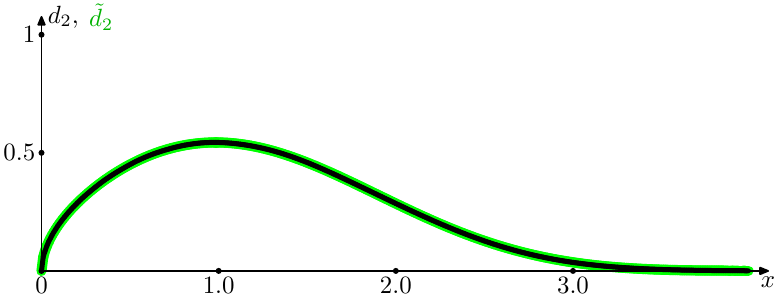}
		\caption{$p=1/3$, $q=1/2$}
		\label{fig1b}
	\end{subfigure}
	\hfill
	\begin{subfigure}[b]{0.75\textwidth}
		\includegraphics[width=\textwidth]{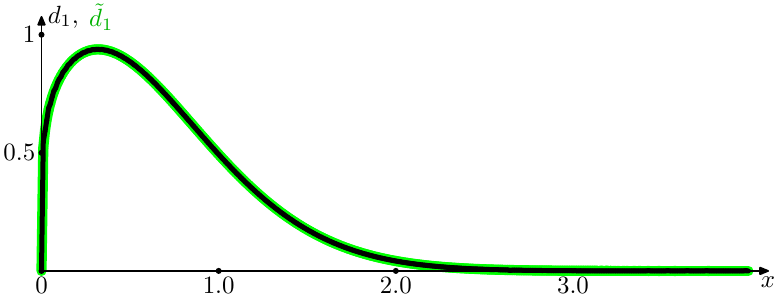}
		\caption{$p=2/3$, $q=1/2$}
		\label{fig1c}
	\end{subfigure}
         \hfill
	\begin{subfigure}[b]{0.75\textwidth}
		\includegraphics[width=\textwidth]{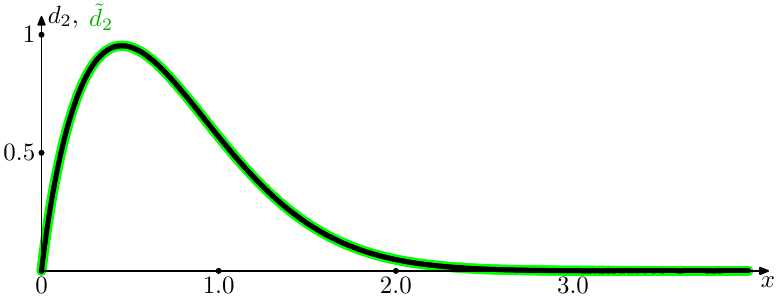}
		\caption{$p=2/3$, $q=1/2$}
		\label{fig1d}
	\end{subfigure}
	\caption{For the case \er{200}, the density \er{010} and its approximation \er{03} are computed.}\label{fig2}
\end{figure}

The results of computations are presented in Fig. \ref{fig2}. For the approximation \er{03}, we use $40\ts40$ of both entries of Taylor coefficients $\pmb{\psi}_{\bf m}$. The function ${\bf K}$, see \er{020}, is computed with the help of the schemes \er{045}, \er{047}, and \er{048}, \er{050}, where the number of iterations $t=250$. For the exact density \er{010}, we use the trapezoidal rule dividing the interval $[0,400\mathbf{i}]$ on $10^5$ nodal points uniformly. Here, we take only $t=50$ iterations in the scheme based on \er{045} and \er{047}. The reason is the large number of nodal points, where $\pmb{\Pi}(z)$ should be computed.
Approximate and ``exact'' plots are almost identical, even if we take $t=50$ iterations for the approximation, instead of $t=250$. We took many iterations just to emphasize that the approximation is calculated very quickly, practically on the fly.

{\section{Conclusion}\lb{sec3}} 
The complete ``power series expansion'' along with the good approximation are derived for the density of {\it martingale limit} of Galton-Watson processes in multidimensional cases. Possible avenues for further research include relaxing conditions {\bf (A)}-{\bf (D)}, with the special focus on {\bf (D)}. If this happens, then we will compute the series for other interesting cases like (c) and (d) from Fig. \ref{fig1}.

\section*{Ethical Statements}

The author has no relevant financial or non-financial interests to disclose. Only standard support for general research work without special interests is listed in the Acknowledgements section below.

The author has no competing interests to declare that are relevant to the content of this article.

The author certifies that he has no affiliations with or involvement in any organization or entity with any financial interest or non-financial interest in the subject matter or materials discussed in this manuscript.

The author has no financial or proprietary interests in any material discussed in this article.

\section*{Data availability statement}

Data sharing not applicable to this article as no datasets were generated or analyzed during the current study. The program code for procedures and functions considered in the article is available from the corresponding author on reasonable request.

\section*{Acknowledgements} 
This paper is a contribution to the project M3 of the Collaborative Research Centre TRR 181 "Energy Transfer in Atmosphere and Ocean" funded by the Deutsche Forschungsgemeinschaft (DFG, German Research Foundation) - Projektnummer 274762653. 
I would like to thank Nick Bingham for discussing some results regarding the multi-dimensional case.

\bibliographystyle{abbrv}

\end{document}